\documentclass[10pt]{article}
\usepackage[usenames]{color} %used for font color
\usepackage{amssymb} %maths
\usepackage{amsmath,graphicx} %maths
\usepackage[utf8]{inputenc} %useful to type directly diacritic characters
\usepackage{hyperref}
\usepackage{cleveref}
\def\R{\mathbb R}

\title{Analytical travelling vortex  solutions of hyperbolic equations  for validating very high order schemes} 

\author{Mario Ricchiuto\thanks{Team CARDAMOM, Inria Bordeaux Sud-Ouest, - 200 av.  de la vieille tour, 33405 Talence, France}\, and Davide Torlo\thanks{Team CARDAMOM, Inria Bordeaux Sud-Ouest, - 200 av.  de la vieille tour, 33405 Talence, France}}

\date{\today}

\begin{document}
\maketitle

\begin{abstract}
Testing the order of accuracy of (very) high order methods for shallow water (and Euler) equations is a delicate operation and the test cases are the crucial starting point of this operation. We provide a short derivation of vortex-like analytical solutions in 2 dimensions for the shallow water equations (and, hence, Euler equations) that can be used to test the order of accuracy of numerical methods. These solutions have different smoothness in their derivatives (up to $\mathcal C^\infty$) and can be used accordingly to the order of accuracy of the scheme to test.
\end{abstract}

\section{Moving vortex   solutions and regularity requirements}

\subsection{Shallow water equations and moving profiles} 

We consider the shallow water equations (SWEs) on a flat bathymetry reading
\begin{equation}\label{eq:00}
\begin{cases}
\partial_t h  & +\nabla\cdot(h \vec u) =0, \\
\partial_t (h \vec u )  & + \nabla\cdot(h\vec u\otimes\vec u)   + gh \nabla h=0,
\end{cases}
\end{equation}
$h : \R^2 \to \R^+$, $\vec u :\R^2 \to \R^2$, $g\in \R^+$.

\noindent
Following \cite{ricchiuto2009stabilized} we consider solutions  of the form $h=H_0(\vec\zeta)$ and $\vec u=\vec u_{\infty}+\vec U_0(\vec\zeta)$,  with 
$\zeta =\vec x - \vec u_{\infty} t$, and $\vec u_{\infty}$ constant.  Replacing in the SWEs we obtain
\begin{equation}\label{eq:00a}
\begin{split}
\begin{cases} \vec U_0 \cdot\nabla_{\vec \zeta\,} H_0   &  + H_0\nabla_{\vec\zeta\,}\cdot\vec U_0 =0,\\
\vec U_0 \cdot\nabla_{\vec \zeta\,} \vec U_0  &   +  \nabla_{\vec\zeta\,}H_0 =0,\\
\end{cases}
\end{split}
\end{equation}
which justifies looking for stationary solutions with solenoidal velocity fields with depth variations uniquely in the cross-stream direction.
These conditions are easily met for solutions with cylindrical symmetry.

\subsection{SWEs  in cylindrical coordinates}
We consider cylindrical coordinates defined in 2D by the distance from the origin $r^2=x^2+y^2$, and the counter-clockwise angle $\theta$
measured from the positive $x$ axis   so that  
\begin{equation}\label{eq:0a}
x = r \cos\theta\;,\;\;y = r \sin\theta.
\end{equation}
We also introduce the direction vectors $\hat r=(\cos\theta,\,\sin\theta)$, and $\hat r^{\perp} =(- \sin\theta,\, \cos\theta)$. 
This notation can be used to write the shallow water equations in polar coordinates as  
\begin{equation}\label{eq:1}
\begin{cases}
\partial_t h& +  \partial_r(h u_r) +\dfrac{1}{r}\left(  \partial_{\theta}(h u_{\theta})  + h u_r\right) =0 \\
\partial_t (h u_r)& +  \partial_r(h u_r^2)   + gh \partial_r h +\dfrac{1}{r}\left(  \partial_{\theta}(h u_r u_{\theta})  + h (u_r^2-u_{\theta}^2)\right) =0 \\
\partial_t (h u_{\theta})& +  \partial_r(h u_ru_{\theta})   +\dfrac{1}{r}\left(   gh \partial_{\theta} h  + \partial_{\theta}(h  u_{\theta}^2)  +2 h u_r u_{\theta} \right) =0 
\end{cases}
\end{equation}
where
\begin{equation}\label{eq:0b}\begin{split}
		&u_r= \vec u\cdot \hat r = \cos(\theta)u_x+\sin(\theta)u_y \,,\;\;
	u_{\theta}= \vec u\cdot \hat r^{\perp} = -\sin(\theta)u_x+\cos(\theta)u_y,\\
	&u_x= \cos(\theta)u_r-\sin(\theta)u_\theta \,,\;\;u_y= \sin(\theta)u_r+\cos(\theta)u_\theta.
	\end{split}
\end{equation}

\subsection{Stationary vortex ODE} 
To mimic \eqref{eq:00a}, we consider the particular case with all zero time derivatives and only radial velocity, i.e., 
\begin{equation}\label{eq:2}
\begin{split}
\partial_t h  =& 0,\\  \partial_t u_r=&\partial_t u_{\theta}=0, \\
h =&h(r)\;,\;\; \partial_{\theta}h =0, \\
u_{\theta} =&u_{\theta}(r)\;,\;\; \partial_{\theta}u_{\theta} =0, \\
u_r =&0.
\end{split}
\end{equation}

With the above hypotheses we can readily check that the first and the last in \eqref{eq:1} are identically  satisfied, while the second reduces to
\begin{equation}\label{eq:3}
h'(r)= \dfrac{u_{\theta}^2}{gr}.
\end{equation}

Assuming further that  $u_{\theta}=\omega(r) r$ we end up with

\begin{equation}\label{eq:4}
h'(r)= \dfrac{r\omega^2(r)}{g}.
\end{equation}
Given a law for the angular velocity $\omega$ this ODE can be integrated to obtain closed form expressions for the depth and, conversely, given a law for $h$, $u_\theta$ and $\omega$ can be obtained differentiating $h$.

\subsection{Regularity requirements for validating high order methods}

To validate higher order methods, exact solutions of the shallow water system should have  enough regularity to allow the  validity of high order approximation results.
The classical interpolation estimate for finite element approximations \cite[\S1.5]{ern2013theory} for a function with $\partial_{\alpha}v\in L^{p}(\Omega)$, for a multi-index $|\alpha|\le l+1$,
is the following

\begin{equation}\label{eq:fe0}
\|v - v_h\|_{L^{p}(\Omega)} \le c h^{l+1} |v|_{l+1,L^{p}(\Omega)},
\end{equation}

with $h$ the mesh size. This means that to benchmark an $l+1$ order accurate method, we need an exact solution with integrable $l+1$ derivatives. In practice, for a system different variables may have
different regularity. This may lead in practice to convergence rates somewhat in between those of the different variables, depending on how the error is defined and measured.

\subsection{Extension to Euler equations}
Similarly to SWEs other moving vortexes can solve exactly Euler equations.
They read
\begin{equation}
	\begin{cases}
		\partial_t \rho &+  \partial_x(\rho u_x) + \partial_y (\rho u_y)=0,\\
		\partial_t (\rho u_x) &+ \partial_x (\rho u_x^2 + p) + \partial_y (\rho u_x u_y)=0,\\
		\partial_t (\rho u_y) &+\partial_x (\rho u_xu_y) + \partial_y (\rho u_y^2 + p)=0,\\
		\partial_t (\rho E) &+ \partial_x \left[ u_x(\rho E+ p) \right] + \partial_y \left[ u_y(\rho E+ p) \right]=0,
	\end{cases}
\end{equation}
where the perfect gas equation of state closes the system, i.e.,
\begin{equation}
	\rho E = \frac{p}{\gamma -1}+\frac{1}{2}\rho (u_x^2+u_y^2).
\end{equation}
Here $\gamma \in \R^+$ is the adiabatic constant.
Equivalently, we can write Euler equations in polar coordinates
\begin{equation}
	\begin{cases}
		\partial_t \rho &+  \frac{1}{r}\partial_r(r\rho u_r) + \frac{1}{r}\partial_\theta (\rho u_\theta)=0,\\
		\partial_t (\rho u_r) &+ \frac{1}{r}\partial_r (r(\rho u_r^2 + p)) + \frac{1}{r}\partial_\theta (\rho u_r u_\theta)-\frac{1}{r}(\rho u_\theta^2+p)=0,\\
		\partial_t (\rho u_\theta) &+\frac{1}{r}\partial_r (r\rho u_ru_\theta) + \frac{1}{r}\partial_\theta (\rho u_\theta^2 + p)+\frac{1}{r}\rho u_\theta u_r=0,\\
		\partial_t (\rho E) &+ \frac{1}{r}\partial_r \left[r u_r(\rho E+ p) \right] + \frac{1}{r}\partial_\theta \left[ u_\theta(\rho E+ p) \right]=0.
	\end{cases}
\end{equation}

Using the vortex form \eqref{eq:2}, we have that $u_r$ is set to 0 and $\partial_\theta$ is equal to 0 for all the unknowns. Hence, the system reduces for steady vortexes to
\begin{equation}\label{eq:ruleEuler}
	r \partial_r p = \rho u_\theta^2.
\end{equation}
This form can be easily solved in two situations. 
\subsubsection*{Isentropic case}
The isentropic case ($S=p/\rho^{\gamma}=$constant) for $\gamma >1$ leads to
\begin{equation}
	\frac{\gamma}{\gamma -1} \partial_r (\rho ^{\gamma-1}) = \frac{u_\theta^2}{r}.
\end{equation}

Following the approach of \eqref{eq:3}, if we define any angular velocity $\omega(r)$ and consequentially 
\begin{equation}
	u_r=0, \qquad u_\theta(r):=r\omega(r), \qquad \rho(r):=\left( \rho_0 + \int \frac{\gamma-1}{\gamma} r \omega^2(r) dr \right)^{\frac{1}{\gamma-1}}, \qquad p=\rho^\gamma,
\end{equation}  
where with the integral is meant up to a constant, that can be set with $\rho$ at $r=\infty$.
Clearly, this formulation is equivalent to the SWE ones setting $\gamma=2$ and $g=2$, and all the following derivation can be used in Euler equations.

\subsubsection*{Isochoric vortex}
The second option to obtain a steady vortex is to set a constant density $\rho(r)=\rho_0$, and \eqref{eq:ruleEuler} becomes a differential equation for $p$.
Again, we obtain  
\begin{equation}
	u_r=0, \qquad u_\theta(r):=r\omega(r), \qquad \rho(r)=\rho_0, \qquad p(r):= \int \rho_0 r \omega^2(r) dr.
\end{equation}

\section{The vortex solution of \cite{ricchiuto2009stabilized}}

The following  is an example  often used in literature obtained by setting 
\begin{equation}\label{eq:5}
\omega(r) =\begin{cases}
	\Gamma (1+ \cos(\pi\dfrac{ r}{r_0})) \quad&\text{if } r\le r_0, \\[5pt]
	0\quad&\text{otherwise}.
\end{cases}
\end{equation}
This expression can be fed into \eqref{eq:4} and integrated backwards from $r=r_0$ to $r$ with initial condition $h=h_0$ to obtain
\begin{equation}\label{eq:6}
h=h_0- \dfrac{\Gamma^2r_0^2}{g \pi^2}
\begin{cases}
H(\pi) - H(\pi\dfrac{ r}{r_0})   \quad&\text{if } r\le r_0, \\[5pt]
0\quad&\text{otherwise},
\end{cases} 
\end{equation}
with
\begin{equation}\label{eq:7}
  H(x)= 2\cos x + 2x\sin x+\dfrac{\cos(2x)}{8}+\dfrac{x\sin(2x)}{4} +\dfrac{12x^2}{16}.
\end{equation}

Note that  

\begin{equation}\label{eq:8}
\begin{split}
  H'(\pi)=0\;,\;\; &u_{\theta}'(\pi)=  0,\\ 
  H^{(2)}(\pi)=  0\;,\;\; &u_{\theta}^{(2)}(\pi)= \pi \ne 0,\\ 
    H^{(3)}(\pi)=   0\;,\;\; &u_{\theta}^{(3)}(\pi)\ne 0, \\ 
        H^{(4)} (\pi)=  0 \;,\;\; &u_{\theta}^{(4)}(\pi)\ne 0 ,\\ 
                H^{(5)} (x)=  6\pi\ne 0\;,\;\; &u_{\theta}^{(5)}(\pi)\ne 0 .
  \end{split}
\end{equation}
According to \eqref{eq:fe0}   this solution should allow to obtain at most second order of accuracy, due to the  limited regularity of the velocity.
In practice, some computations have shown a  little extra convergence for the depth, for which, at least on relatively coarse resolutions, one can manage to 
obtain toughly  third order convergence. The convergence study obtained with a WENO5 on Cartesian grids on $[0,1]^2$ with a RK(6,5) is shown in \cref{tab:convergenceCos}. The tests are run on the vortex defined with $r_0=0.25$ and centered in $(0.5,0.5)$, $h_0=1$ and $\Gamma$ such that $h(0)=0.99$. It is clear that the order of accuracy for $h$ cannot reach 3 and that the velocities can be approximated with only order 2.

\begin{table}
	\centering
\begin{tabular}{|c||c|c||c|c||c|c||}
	\hline
	  N   & Error h  &  Order h & Error u  &  Order u & Error v  &  Order v \\
    8  &   2.755e-04  &  0.000  &  3.072e-03 & 0.000  &  3.072e-03 & 0.000 \\ 
   16  &   1.650e-04  &  0.739  &  8.940e-04 & 1.781  &  8.938e-04 & 1.781 \\ 
   32  &   3.039e-05  &  2.441  &  1.654e-04 & 2.434  &  1.654e-04 & 2.434 \\ 
   64  &   4.188e-06  &  2.859  &  3.568e-05 & 2.213  &  3.568e-05 & 2.213 \\ 
  128  &   5.018e-07  &  3.061  &  7.879e-06 & 2.179  &  7.879e-06 & 2.179 \\ 
  256  &   5.755e-08  &  3.124  &  1.588e-06 & 2.311  &  1.588e-06 & 2.311 \\ 
  512  &   6.320e-09  &  3.187  &  2.970e-07 & 2.418  &  2.970e-07 & 2.418 \\ 

	\hline
\end{tabular}
\caption{Order of convergence for \eqref{eq:7} obtained with WENO5 on a cartesian grid with $N\times N$ cells\label{tab:convergenceCos}}
\end{table}

\section{Non compact supported vortexes}
There are many other vortexes for Euler's equation, Euler himself in \cite{euler1775principes} suggested a similar simplification to obtain a steady solution, and in 1998 Shu \cite{shu1998essentially} proposed a vortex to study the accuracy of WENO5 schemes from a qualitatively point of view. 
In the following years, this tests and some modifications of it became a real benchmark for testing the accuracy of many schemes \cite{hesthaven2007nodal,wang2007spectral,wang2009unifying,abgrall2019high,abgrall2020high}.
The main idea follows form the definition of a vortex where the base function $\omega(r)$ is a Gaussian function $\omega(r) = \Gamma e^{-(r/r_0)^2}$, where $r_0$ is a rescaling factor for the width and $\Gamma$ rescales the amplitude.
Despite being the benchmark test for many problems, this test does not have a compact support and boundaries are a real issue, in particular when dealing with very high order of accuracy methods, as shown, for instance, in \cite{spiegel2015survey}. Nevertheless, it was often used with periodic boundary conditions even with nonzero background speed. This leads to traveling discontinuities (on the derivatives) from the boundaries to all over the domain.
The techniques to \textit{cure} this issue were various starting from enlarging the domain (or equivalently reducing $r_0$) to treating with nonperiodic boundary conditions (inflow/outflow with steady vortex). Of course increasing the domain means that the computational costs increase too and using inflow/outflow does not allow to make the vortex travel outside the domain.

We can easily compute 
\begin{equation}\label{eq:Gauss}
	\begin{cases}
		\omega(r)= \Gamma e^{-(r/r_0)^2}\\
		u_\theta(r) = r \omega(r)\\
		h(r) = h_0-\frac{\Gamma^2 r_0^2}{4g}e^{-2(r/r_0)^2}=h_0-\frac{r_0^2}{4g}\omega^2(r).
	\end{cases}
\end{equation}

\begin{figure}
	\centering
	\includegraphics[height=0.95\textheight,trim={50 70 50 50},clip]{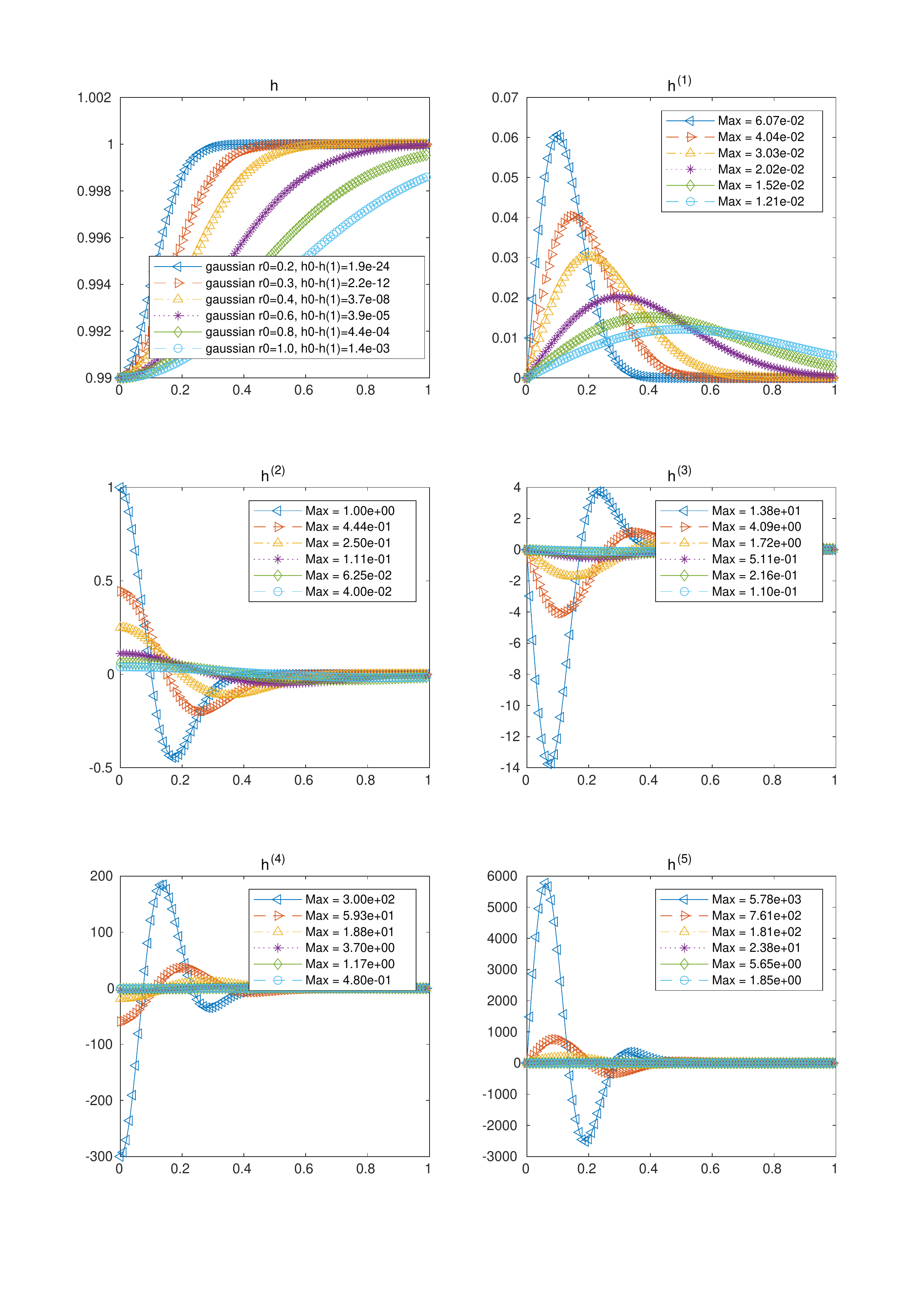}
	\caption{$h$ profile and its derivatives for vortexes \eqref{eq:Gauss} with $r_0=0.2,0.3,0.4,0.6,0.8,1$\label{fig:hGauss}}
\end{figure}

\begin{figure}
	\centering
	\includegraphics[height=0.95\textheight,trim={50 70 50 50},clip]{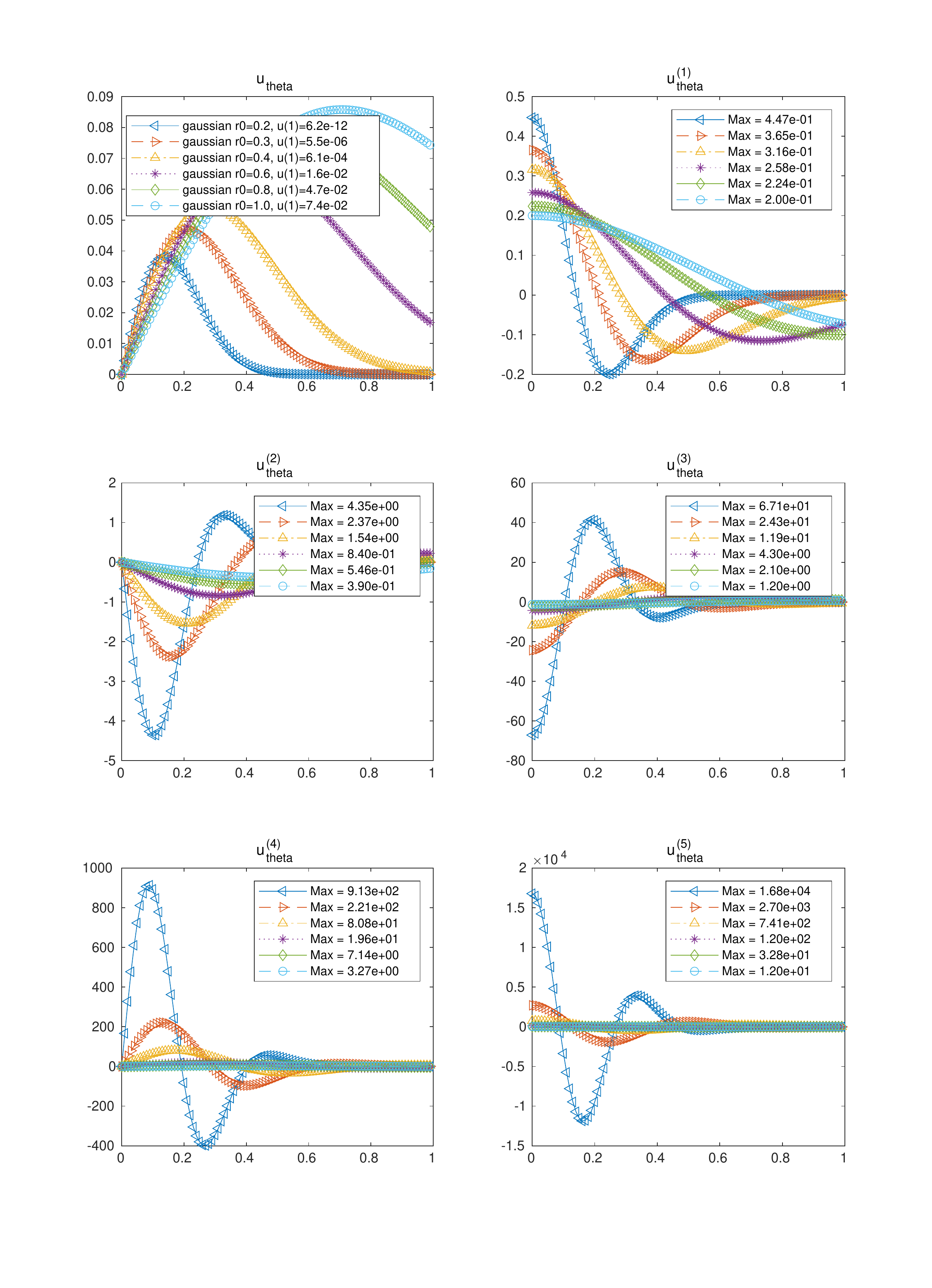}
	\caption{$u_\theta$ profile and its derivatives for vortexes \eqref{eq:Gauss} with $r_0=0.2,0.3,0.4,0.6,0.8,1$\label{fig:uGauss}}
\end{figure}

We plot the profile of $h(r)$ and $u_\theta(r)$ for $h_0=1$, $\Gamma$ such that $h(0)=0.99$ and different $r_0$ in \cref{fig:hGauss,fig:uGauss}.

On the pictures one finds also the values of the maximum absolute values of the derivatives and the values of speed and height at $r=1$. As one can see, increasing $r_0$ the derivatives decrease, but the values at $r=1$ are diverging from $h_0$ for $h$ and 0 for $u_\theta$. This means that boundary effects might disturb the convergence process.
For the numerical part, we will test only $r_0=0.2,\,0.3,\,0.4$, as all the other parameters leads to either very high derivatives or boundary effects. This is also visible in these tests, but this is the range where they are more bounded. 

\section{Compact supported traveling vortexes of arbitrary smoothness}
\subsection{Iterative correction of the RB-vortex to obtain  arbitrary smoothness} 
In order to obtain more regularity in the vortex, and to be able to test higher order accuracy methods, we can generalize the previous solution. We remark that  for $r/r_0\le 1$ \eqref{eq:5} is equivalent to
\begin{equation}\label{eq:9}
\omega = \Gamma \left(1+\cos(\rho) \right)= 2 \Gamma  \cos^2{(\rho/2)}\;,\;\;\rho =\pi\dfrac{ r}{r_0}.
\end{equation}
A natural way to improve the regularity of this definition for $\rho=\pi$ is to increase the exponent of the cosinus.
So we look into definitions of the type 
\begin{equation}\label{eq:10}
\omega = 2^p \Gamma  \cos^{2p}{(\rho/2)} = \Gamma (1+ \cos{\rho})^p\;,\;\; p\ge 1 
\end{equation}
which allows to increase the regularity, keeping bounded values of the $(p+2)$th derivative.\\

In practice, we need to integrate  the ODE
\begin{equation}\label{eq:11}
h'(r) = \dfrac{ 4^p \Gamma^2  }{g} r\,\cos^{4 p}(\pi\dfrac{ r}{2 r_0}  )\;,\;\; r\in[0,\,r_0],
\end{equation}
with the condition $h(r_0)=h_0$. The solution to this problem within $[0,\,r_0]$ can be written as 
\begin{equation}\label{eq:12}
h(r) = h_0 - \dfrac{1}{g}\left(\dfrac{2^p\Gamma r_0}{\pi}\right)^2  (H_{p}(\pi/2) - H_p(\rho/2)) \;,\;\;
\rho = \pi\dfrac{ r}{r_0} ,
\end{equation}
having set  
\begin{equation}\label{eq:13}
H_p(x) =   \int^x y \cos^{4p}(y) \,dy
\end{equation}
up to an additional constant. Comparing to \eqref{eq:6} we can deduce that
\begin{equation}\label{eq:14}
%h=h_0- \dfrac{ (4\Gamma r_0)^2}{g \pi^2} \left( \dfrac{H(2  \pi/2) - H(2 \rho/2) }{16} \right)\;,\;\;
%H(x)= 2\cos x + 2x\sin x+\dfrac{\cos(2x)}{8}+\dfrac{x\sin(2x)}{4} +\dfrac{12x^2}{16}
H_1(x) =\dfrac{1}{16}H(2 x)
\end{equation}
with $H(x)$  defined in \eqref{eq:7}. For $p>1$ we can  use the iterative integration rule
\begin{equation}\label{eq:15}
 \int \cos^n(x)\, dx = \dfrac{\cos^{n-1}(x)\sin (x)}{n} +\dfrac{n-1}{n} \int \cos^{n-2}(x)\,dx
 \end{equation}
 to show that 
 \begin{equation*}%\label{eq:16}
 \begin{split}
H_p(x):=&\int x \cos^{4p}(x)=x\int  \cos^{4p}(x) -\int\int  \cos^{4p}(x) \\
=& x\dfrac{\cos^{4p-1}(x)\sin(x) }{4p} +x\dfrac{4p-1}{4p}\int  \cos^{4p-2}(x)
-\int\left\lbrace \dfrac{\cos^{4p-1}(x)\sin(x) }{4p} +\dfrac{4p-1}{4p}\int  \cos^{4p-2}(x)
 \right\rbrace\\
= &x\dfrac{\cos^{4p-1}(x)\sin(x) }{4p} - \int\dfrac{\cos^{4p-1}(x)\sin(x) }{4p}
+\dfrac{4p-1}{4p} \left\{
x \int  \cos^{4p-2}(x) -  \int \int \cos^{4p-2}(x) 
\right\}\\
%=& x\dfrac{\cos^{4p-1}(x)\sin(x) }{4p} + \dfrac{\cos^{4p}(x)}{(4p)^2}%\dfrac{1}{4p}\int u^{4p-1} du 
%+  \dfrac{4p-1}{4p}\int x \cos^{4p-2}(x)\\
=& x\dfrac{\cos^{4p-1}(x)\sin(x) }{4p} + \dfrac{\cos^{4p}(x)}{(4p)^2}%\dfrac{1}{4p}\int u^{4p-1} du 
+ \\
&
 \dfrac{4p-1}{4p}\left\{
x\dfrac{\cos^{4p-3}(x)\sin(x) }{4p-2}+ \cos^{4p-2}(x)
+\dfrac{4p-3}{4p-2}\int x \cos^{4p-4}(x)
\right\}\\
=& 
\dfrac{4p-1}{4p}\dfrac{4p-3}{4p-2}H_{p-1}(x)
+ x\dfrac{\cos^{4p-3}(x)\sin(x) }{4p}
\left\{
\dfrac{4p-1}{4p-2} + \cos^2(x)
\right\}
%+ x\dfrac{\cos^{4p-1}(x)\sin(x) }{4p}
%+  \dfrac{4p-1}{4p}x\dfrac{\cos^{4p-3}(x)\sin(x) }{4p-2}
+\\
&  \cos^{4p-2}(x)\left\{
 \dfrac{\cos^{2}(x)}{(4p)^2}+  \dfrac{4p-1}{4p(4p-2)^2} \right\}.
%=& x\dfrac{\cos^{4p-1}(x)\sin(x) }{4p} +   \cos^{4p}(x)
%+\dfrac{4p-1}{4p}\,\mathcal{I}_{4p-2}\\
%%%%
%%%%
%%%%
%=& x\dfrac{\cos^{4p-3}(x)\sin(x)}{4p} \left\{   \dfrac{4p }{4p-2} +\cos^2(x)   \right\}
%+   \cos^{4p-2}(x) (1+ \cos^2(x) )
%+\dfrac{4p-1}{4p}  
%\dfrac{4p-3}{4p-2}  H_{p-1}
\end{split}
 \end{equation*}
This formula can be used to compute the exact solution  iteratively for any $p\ge1$ with initial value given by \eqref{eq:14}.
In \cref{app:exactCos}, we provide some values of $H_p$.

\begin{figure}
	\centering
	\includegraphics[height=0.95\textheight,trim={50 70 50 50},clip]{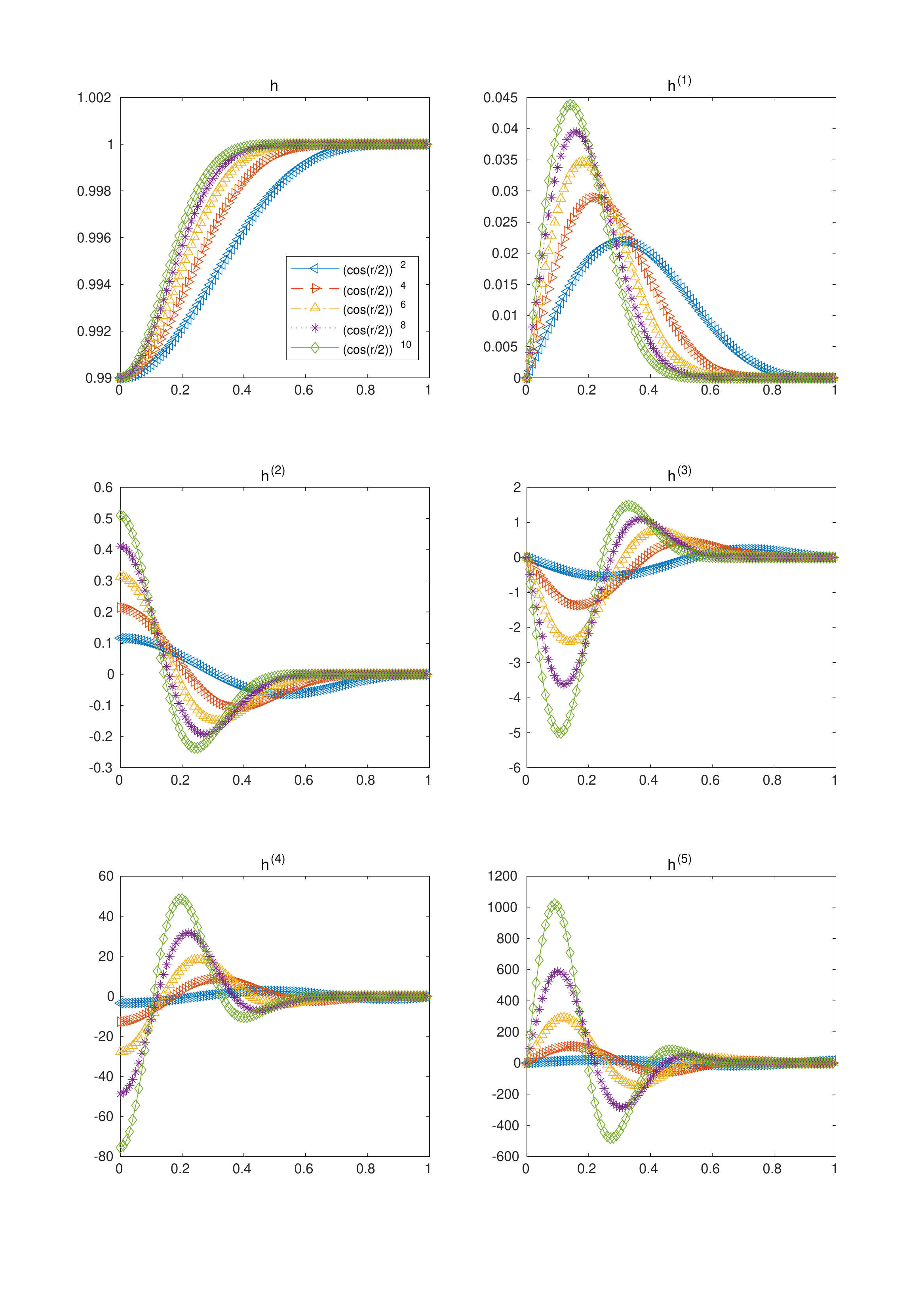}
	\caption{$h$ profile and its derivatives for vortexes \eqref{eq:12} with $p=1,\dots , 5$\label{fig:hCos}}
\end{figure}

\begin{figure}
	\centering
	\includegraphics[height=0.95\textheight,trim={50 70 50 50},clip]{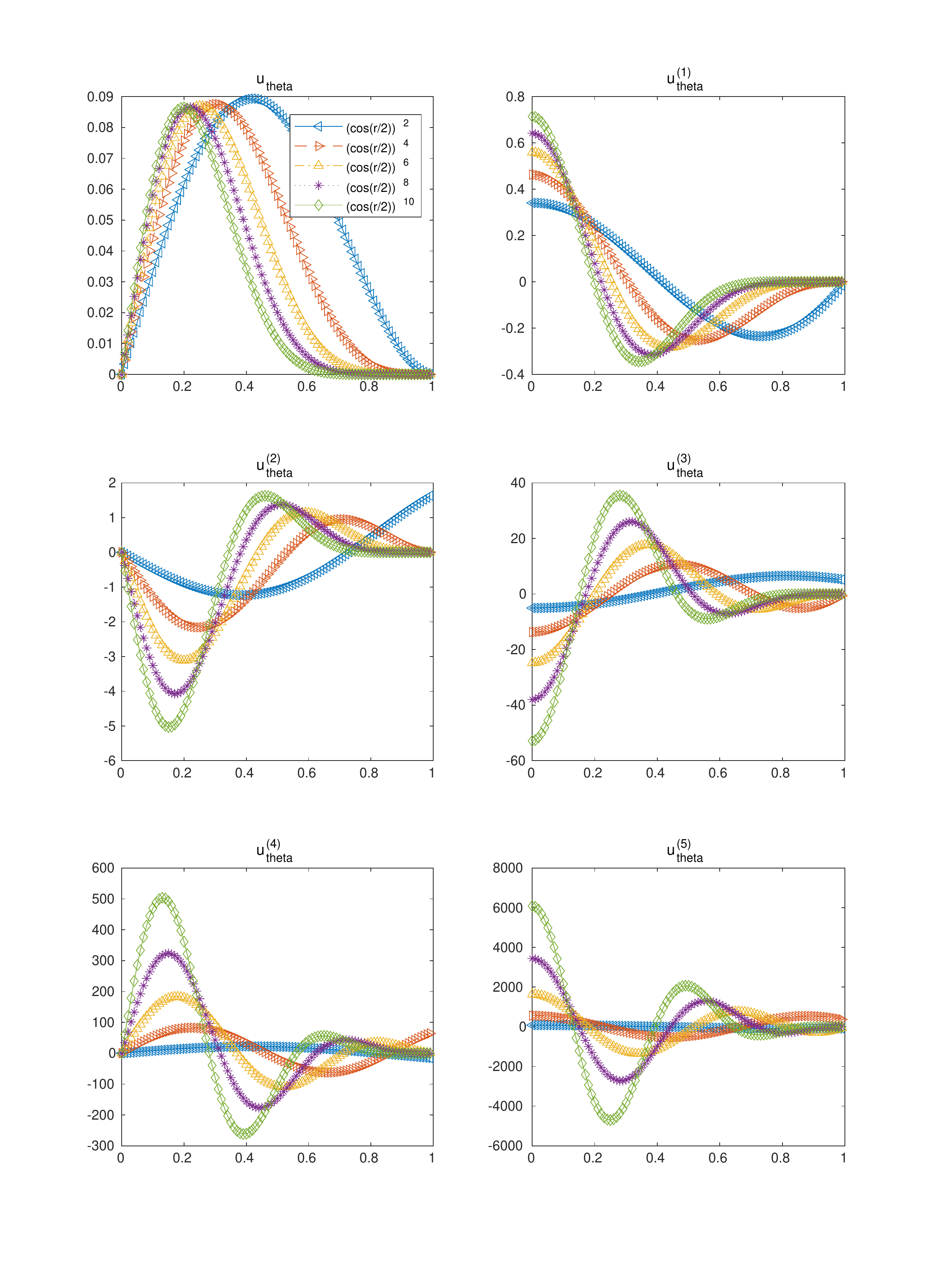}
	\caption{$u_\theta$ profile and its derivatives for vortexes \eqref{eq:12} with $p=1,\dots , 5$\label{fig:uCos}}
\end{figure}

These vortexes suited for testing high order methods, as they are $\mathcal C^{2p}$. Nevertheless, it is not recommended to use a very high $p$ for not so high order methods. Indeed, we see that the estimation on the error \eqref{eq:fe0} depends not only on the smoothness of the solution, but also on the seminorms. Essentially, if the derivatives are too large, we might need very fine meshes to reach the expected order of accuracy.
To assess this information we plot for various vortexes the profile of $h(r)$ and of $u_\theta(r)$ and their derivatives (up to the 5th derivative). We fix for all the vortexes $r_0=1$, the values $h_0=1$ and $h(0)=0.99$, by setting $\Gamma$ and $g=1$. 
In \cref{fig:hCos,fig:uCos} the plot related to the vortexes of type \eqref{eq:12} for variables $h$ and $u_\theta$ respectively.

In \cref{fig:hCos} we barely see that the fifth derivative of the case $p=1$ is not zero at $r=1$, but we observe that the amplitude of the derivative functions increases with $p$, in particular the fifth derivative of $p=5$ is 10 times larger than $p=2$.
In \cref{fig:uCos} we immediately see that the second derivative of $u_\theta$ for $p=1$ is discontinuous in $r=1$, and the fourth derivative is discontinuous also for $p=2$. Again, the higher we choose $p$ the larger the derivatives become.

\subsection{A few $C^{\infty}$ examples} 
In this case we start from \eqref{eq:4} and we reverse it to obtain a definition of the angular velocity  given the depth:
\begin{equation}
\label{eq:d0}
	\omega = \sqrt{\frac{g h' }{r}}.
\end{equation}
To avoid the singularity in $r=0$, we define the depth as a function of $\rho := (r/r_0)^2$. An example of
  $C^{\infty}$ compactly supported function is  obtained with the definition
\begin{equation}
	h=h_0-\Gamma^2\begin{cases}
		 e^{-\dfrac{1}{(1-\rho)^2}} & \text{ if }\rho < 1,\\[10pt]
		0 & \text{ else},
	\end{cases}\qquad \rho =\left(\dfrac{r}{r_0}\right)^2.
\end{equation}
%$h$ and $\rho$ would be then $\mathbb C^\infty((0,\infty))$ and hence, no problem appears on the circumference $r=r_0$.
The angular velocity, and thus the tangential linear velocity, can be readily computed from \eqref{eq:d0}:
 \begin{equation}
	\omega = 
%	\begin{cases}
%		2 \sqrt{\frac{C g e^{-\frac{1}{(1-(r/r_0)^2)^2}} }{r_0 \left( 1-(r/r_0)^2\right)^{3} }}  & \text{ if }r<r_0\\
%		0 & \text{ else.}	
%	\end{cases}=
 \begin{cases}
	 2\Gamma e^{-\dfrac{1}{2(1-\rho)^2}}\sqrt{\dfrac{g}{r_0\left( 1-\rho\right)^{3} }}    & \text{ if }r<r_0,\\[10pt]
	0 & \text{ else},	
\end{cases}\qquad \rho =\left(\dfrac{r}{r_0}\right)^2.
\end{equation}

Though being $\mathcal C^\infty$, the previous function has very large derivatives and need very fine mesh to see the order of the mesh.

Alternatives could be obtained, for instance, using higher power in the exponential, for example
\begin{equation}\label{eq:hvortexCInfp}
	h=h_0-\Gamma^2\begin{cases}
		e^{-\dfrac{1}{(1-\rho)^p}} & \text{ if }\rho < 1,\\[10pt]
		0 & \text{ else},
	\end{cases}\qquad \rho =\left(\dfrac{r}{r_0}\right)^2,
\end{equation}
leading to
 \begin{equation}
	\omega = 
	%	\begin{cases}
	%		2 \sqrt{\frac{C g e^{-\frac{1}{(1-(r/r_0)^2)^2}} }{r_0 \left( 1-(r/r_0)^2\right)^{3} }}  & \text{ if }r<r_0\\
	%		0 & \text{ else.}	
	%	\end{cases}=
	\begin{cases}
		\Gamma \sqrt{\dfrac{2pg}{r_0\left( 1-\rho\right)^{p+1} }} e^{-\dfrac{1}{2(1-\rho)^p}}   & \text{ if }r<r_0,\\[10pt]
		0 & \text{ else},	
	\end{cases}\qquad \rho =\left(\dfrac{r}{r_0}\right)^2.
\end{equation}

Or one can make the derivatives a bit smaller with an additional $\arctan$ function, i.e.,
\begin{equation}\label{eq:hvortexCTan}
	h=h_0-\Gamma^2\begin{cases}
		e^{-\dfrac{1}{\arctan^p(1-\rho)}} & \text{ if }\rho < 1,\\[10pt]
		0 & \text{ else},
	\end{cases}\qquad \rho =\left(\dfrac{r}{r_0}\right)^2,
\end{equation}
leading to
\begin{equation}
	\omega = 
	%	\begin{cases}
	%		2 \sqrt{\frac{C g e^{-\frac{1}{(1-(r/r_0)^2)^2}} }{r_0 \left( 1-(r/r_0)^2\right)^{3} }}  & \text{ if }r<r_0\\
	%		0 & \text{ else.}	
	%	\end{cases}=
	\begin{cases}
		\Gamma e^{-\dfrac{1}{2\arctan(1-\rho)^p}}\sqrt{\dfrac{2pg}{r_0\arctan^{p+1}\left( 1-\rho\right) }\frac{1}{1+(1-\rho)^2}}    & \text{ if }r<r_0,\\[10pt]
		0 & \text{ else},	
	\end{cases}\qquad \rho =\left(\dfrac{r}{r_0}\right)^2.
\end{equation}

%\section{Vortex shapes and magnitude of the derivatives}

%The derivatives  of the exact solution play a major role in the error estimate \eqref{eq:fe0}.  Solutions with sharply varying high  order derivatives may require fine meshes
%to resolve the  relevant quantities dominating the error  coefficient before providing the asymptotic convergence sought when validating high order schemes.
%Exact solutions allowing this asymptotic range to settle in sooner are   more interesting as this range is in practice observable on a wider window of  mesh sizes before machine accuracy is  eventually reached.
%The objective of this section is thus to look into the behaviour of the solution derivatives for  the vortices proposed above.
%

\begin{figure}
	\centering
	\includegraphics[height=0.95\textheight,trim={50 70 50 50},clip]{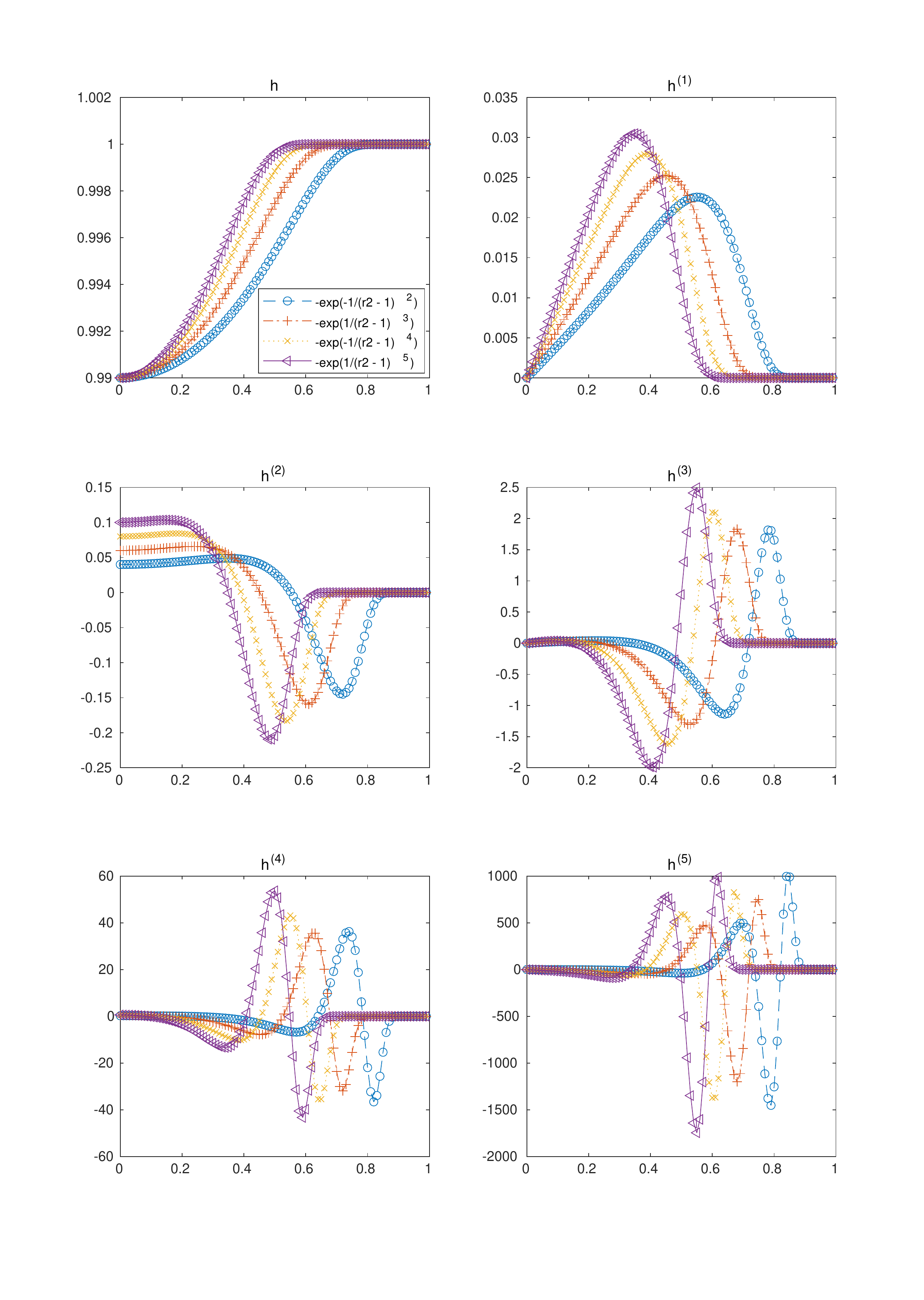}
	\caption{$h$ profile and its derivatives for vortexes \eqref{eq:hvortexCInfp} with $p=2,\dots , 4$\label{fig:hCInf}}
\end{figure}

\begin{figure}
	\centering
	\includegraphics[height=0.95\textheight,trim={50 70 50 50},clip]{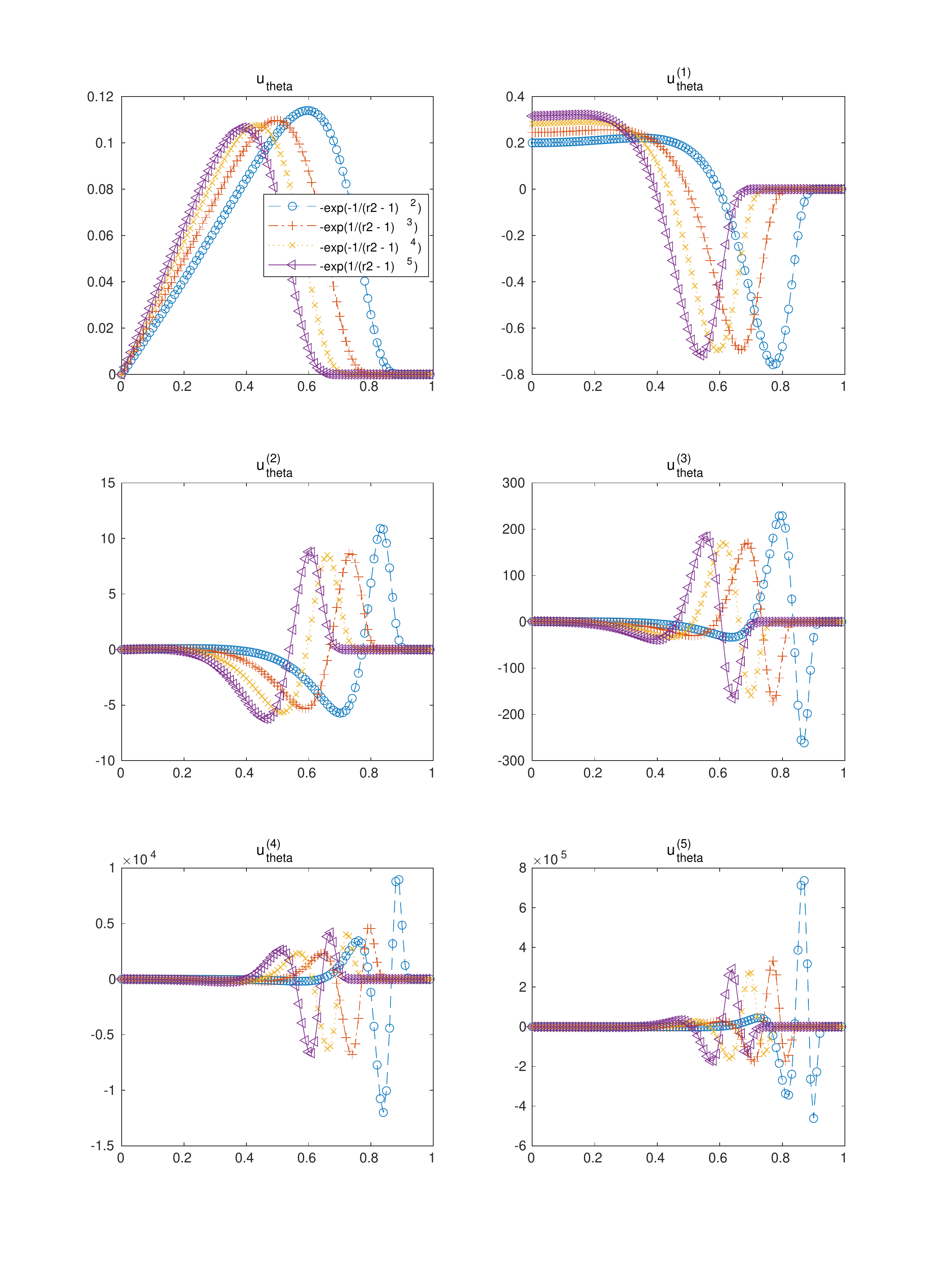}
	\caption{$u_\theta$ profile and its derivatives for vortexes \eqref{eq:hvortexCInfp} with $p=2,\dots , 4$\label{fig:uCInf}}
\end{figure}

\begin{figure}
	\centering
	\includegraphics[height=0.95\textheight,trim={50 70 50 50},clip]{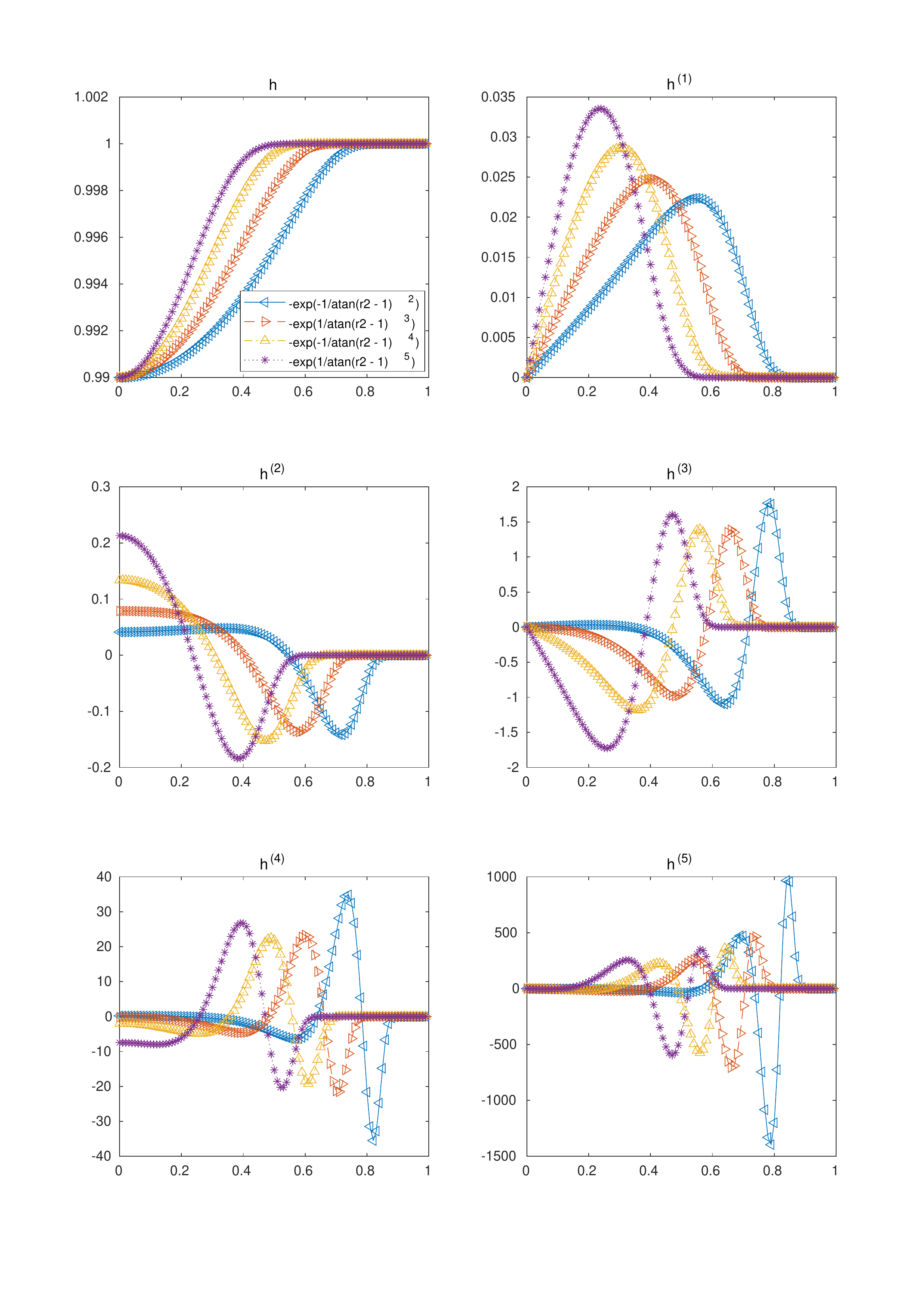}
	\caption{$h$ profile and its derivatives for vortexes \eqref{eq:hvortexCTan} with $p=2,\dots , 4$\label{fig:hCTan}}
\end{figure}

\begin{figure}
	\centering
	\includegraphics[height=0.95\textheight,trim={50 70 50 50},clip]{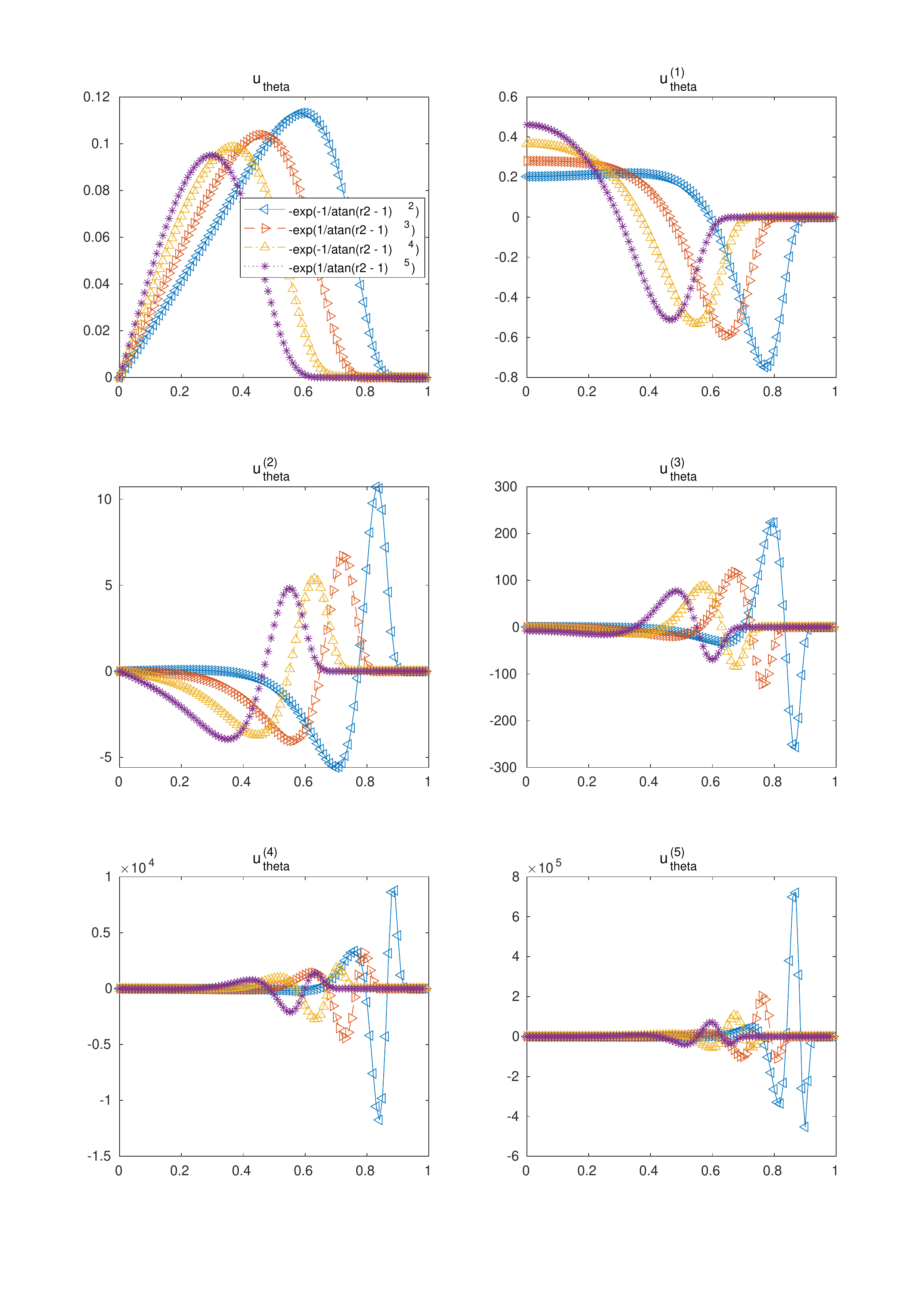}
	\caption{$u_\theta$ profile and its derivatives for vortexes \eqref{eq:hvortexCTan}  with $p=2,\dots , 4$\label{fig:uCTan}}
\end{figure}

As for the previous examples, though being all these vortexes suited for testing high order methods, as they are $\mathcal C^\infty$, the amplitude of their derivatives changes considerably.
To assess this information we plot for various vortexes the profile of $h(r)$ and of $u_\theta(r)$ and their derivatives (up to the 5th derivative). We fix for all the vortexes $r_0=1$, the values $h_0=1$ and $h(0)=0.99$, by setting $\Gamma$ and $g=1$. 
In \cref{fig:hCInf,fig:uCInf} we plot the vortexes of type \eqref{eq:hvortexCInfp} for variables $h$ and $u_\theta$ respectively, while in \cref{fig:hCTan,fig:uCTan} we plot the vortexes for \eqref{eq:hvortexCTan}.

%In \cref{fig:hCos} we barely see that the fifth derivative of the case $p=1$ is not zero at $r=1$, but we observe that the amplitude of the derivative functions increases with $p$, in particular the fifth derivative of $p=5$ is 10 times larger than $p=2$.
%In \cref{fig:uCos} we immediately see that the second derivative of $u_\theta$ for $p=1$ is discontinuous in $r=1$, and the fourth derivative is discontinuous also for $p=2$. Again, the higher we choose $p$ the larger the derivatives become.

With the $\mathcal C^\infty$ vortexes \eqref{eq:hvortexCInfp} and \eqref{eq:hvortexCTan}, the profile of $h$ grows to $h_0$ more slowly, see \cref{fig:hCInf}, but the high derivatives are larger with respect to the ones of the previous examples, but they do not vary much with respect to the parameter $p$. Also for $u_\theta$ we observe similar behaviors: comparing the 5th derivatives of these vortexes with the previous ones, we have values at least 100 times larger.
With the $\arctan$ test cases we actually can slightly decrease the amplitude of high derivatives for large $p$, as shown in \cref{fig:uCTan}.

\section{Numerical tests}
We test a 5th order method with the presented vortexes. The method consists of a finite volume discretization on a Cartesian grid with WENO5 reconstruction and Rusanov numerical flux. We use 4 Gauss-Legendre points for one-dimensional quadrature rules. The time discretization is carried out with Butcher's RK(6,5) 6 stages 5th order method, see \cref{sec:RK65} for its Butcher's tableau.  The CFL number used is 0.95. The domain is $[0,1]^2$ and it is discretized with uniform Cartesian grids and periodic boundary conditions. The vortexes are set with $h_0=1$, $h_{min}=0.99$ and $r_0=0.45$, while for the non compact supported one we test different $r_0$ parameters. Final time is set to $T=1$, but it might be useful to run the code for larger or smaller times. In \cite{spiegel2015survey} it is well observed how this impact on the convergence of the method as the machine precision error accumulates exponentially in time, one might see its effect before reaching the desired accuracy. On the other side, small final time $T$ leads to very small errors also for coarser meshes, not letting appreciate the quality of the high order schemes.

\begin{table}
	\centering
	$r_0$=0.1\\
	\begin{tabular}{|c||c|c||c|c||c|c|}
		\hline
		  N   & Error h  &  Order h & Error u  &  Order u & Error v  &  Order v \\ 
   25  &   1.073e-05  &  0.000  &  3.078e-04 & 0.000  &  3.080e-04 & 0.000 \\ 
   50  &   1.537e-06  &  2.804  &  3.432e-05 & 3.165  &  3.432e-05 & 3.166 \\ 
  100  &   5.777e-08  &  4.734  &  2.174e-06 & 3.980  &  2.174e-06 & 3.980 \\ 
  200  &   2.118e-09  &  4.770  &  1.045e-07 & 4.379  &  1.045e-07 & 4.379 \\ 
  300  &   3.081e-10  &  4.754  &  1.630e-08 & 4.581  &  1.630e-08 & 4.582 \\ 
  400  &   7.620e-11  &  4.857  &  4.354e-09 & 4.589  &  4.354e-09 & 4.589 \\ 
  500  &   2.632e-11  &  4.765  &  1.568e-09 & 4.575  &  1.568e-09 & 4.576 \\ 
  600  &   1.198e-11  &  4.315  &  6.819e-10 & 4.569  &  6.819e-10 & 4.568 \\ 

		\hline
	\end{tabular}\\
	$r_0$=0.15\\
	\begin{tabular}{|c||c|c||c|c||c|c|}
		\hline
		  N   & Error h  &  Order h & Error u  &  Order u & Error v  &  Order v \\ 
   25  &   1.014e-05  &  0.000  &  1.804e-04 & 0.000  &  1.804e-04 & 0.000 \\ 
   50  &   4.589e-07  &  4.466  &  1.336e-05 & 3.755  &  1.336e-05 & 3.755 \\ 
  100  &   1.641e-08  &  4.806  &  7.597e-07 & 4.136  &  7.597e-07 & 4.136 \\ 
  200  &   7.739e-10  &  4.406  &  4.224e-08 & 4.169  &  4.224e-08 & 4.169 \\ 
  300  &   3.228e-10  &  2.156  &  1.373e-08 & 2.772  &  1.373e-08 & 2.772 \\ 
  400  &   2.025e-10  &  1.622  &  8.477e-09 & 1.676  &  8.477e-09 & 1.675 \\ 
  500  &   1.467e-10  &  1.442  &  6.465e-09 & 1.215  &  6.465e-09 & 1.215 \\ 
  600  &   1.127e-10  &  1.447  &  5.399e-09 & 0.988  &  5.399e-09 & 0.988 \\ 

		\hline
	\end{tabular}\\
	$r_0$=0.2\\
	\begin{tabular}{|c||c|c||c|c||c|c|}
		\hline
		  N   & Error h  &  Order h & Error u  &  Order u & Error v  &  Order v \\ 
   25  &   3.422e-06  &  0.000  &  7.886e-05 & 0.000  &  7.886e-05 & 0.000 \\ 
   50  &   2.184e-07  &  3.970  &  9.671e-06 & 3.028  &  9.671e-06 & 3.028 \\ 
  100  &   9.568e-08  &  1.190  &  2.710e-06 & 1.835  &  2.710e-06 & 1.835 \\ 
  200  &   4.185e-08  &  1.193  &  1.467e-06 & 0.885  &  1.467e-06 & 0.885 \\ 
  300  &   2.465e-08  &  1.306  &  1.062e-06 & 0.798  &  1.062e-06 & 0.798 \\ 
  400  &   1.679e-08  &  1.333  &  8.468e-07 & 0.786  &  8.468e-07 & 0.786 \\ 
  500  &   1.241e-08  &  1.356  &  7.097e-07 & 0.791  &  7.097e-07 & 0.791 \\ 
  600  &   9.701e-09  &  1.351  &  6.143e-07 & 0.792  &  6.143e-07 & 0.792 \\ 

		\hline
	\end{tabular}
	\caption{Order of convergence for \eqref{eq:Gauss} with $r_0=0.1, 0.15, 0.2$ obtained with WENO5 on $[0,1]^2$ with a cartesian grid with $N\times N$ cells\label{tab:Gauss}}
\end{table}

With the non compact supported vortex we obtain the results in \cref{tab:Gauss}. We observe for large values of $r_0$ that the error from the boundary effect is predominant already for very coarse meshes. For $r_0=0.2$ (corresponding to $r_0=0.4$ in \cref{fig:hGauss,fig:uGauss}) we have a discontinuity of the order of $10^{-6}$ for the velocity variables on the boundary of the domain (periodic BCs) and this prevent reaching more than 3rd order in the coarse regime. For $r_0=0.15$ we see again boundary effects but for finer meshes, hence it is possible to reach an order of accuracy 4 for few steps in the mesh refinement process.
On the other side, we have for $r_0=0.1$ (corresponding to $r_0=0.2$ in \cref{fig:hGauss,fig:uGauss}) that $h$ at the boundaries is $1$ up to machine precision and $u_\theta\approx 10^{-12}$. Indeed, the convergence in this test is smooth enough for $T=1$, even for not too fine meshes. Nevertheless, one should really be careful using this vortex in hitting machine precision at the boundaries.

\begin{table}
	\centering
	p=1\\
	\begin{tabular}{|c||c|c||c|c||c|c|}
		\hline
		  N   & Error h  &  Order h & Error u  &  Order u & Error v  &  Order v \\ 
   25  &   1.767e-04  &  0.000  &  1.477e-03 & 0.000  &  1.477e-03 & 0.000 \\ 
   50  &   1.716e-05  &  3.365  &  2.420e-04 & 2.610  &  2.420e-04 & 2.610 \\ 
  100  &   1.409e-06  &  3.606  &  4.610e-05 & 2.392  &  4.610e-05 & 2.392 \\ 
  200  &   1.367e-07  &  3.366  &  9.036e-06 & 2.351  &  9.036e-06 & 2.351 \\ 
  300  &   3.817e-08  &  3.147  &  3.416e-06 & 2.399  &  3.416e-06 & 2.399 \\ 
  400  &   1.576e-08  &  3.076  &  1.705e-06 & 2.416  &  1.705e-06 & 2.416 \\ 
  500  &   7.966e-09  &  3.056  &  9.942e-07 & 2.416  &  9.942e-07 & 2.416 \\ 
  600  &   4.568e-09  &  3.050  &  6.415e-07 & 2.403  &  6.415e-07 & 2.403 \\ 

		\hline
	\end{tabular}\\
	p=2\\
	\begin{tabular}{|c||c|c||c|c||c|c|}
		\hline
		  N   & Error h  &  Order h & Error u  &  Order u & Error v  &  Order v \\ 
   25  &   2.072e-04  &  0.000  &  1.615e-03 & 0.000  &  1.614e-03 & 0.000 \\ 
   50  &   2.274e-05  &  3.188  &  1.827e-04 & 3.144  &  1.827e-04 & 3.144 \\ 
  100  &   9.111e-07  &  4.641  &  1.717e-05 & 3.411  &  1.717e-05 & 3.411 \\ 
  200  &   3.305e-08  &  4.785  &  1.136e-06 & 3.918  &  1.136e-06 & 3.918 \\ 
  300  &   4.701e-09  &  4.810  &  2.215e-07 & 4.031  &  2.215e-07 & 4.031 \\ 
  400  &   1.207e-09  &  4.725  &  6.977e-08 & 4.016  &  6.977e-08 & 4.016 \\ 
  500  &   4.413e-10  &  4.511  &  2.870e-08 & 3.981  &  2.870e-08 & 3.981 \\ 
  600  &   2.091e-10  &  4.096  &  1.392e-08 & 3.966  &  1.393e-08 & 3.966 \\ 

		\hline
	\end{tabular}\\
	p=3\\
	\begin{tabular}{|c||c|c||c|c||c|c|}
		\hline
		  N   & Error h  &  Order h & Error u  &  Order u & Error v  &  Order v \\ 
   25  &   1.684e-04  &  0.000  &  1.471e-03 & 0.000  &  1.471e-03 & 0.000 \\ 
   50  &   3.635e-05  &  2.212  &  2.538e-04 & 2.536  &  2.537e-04 & 2.536 \\ 
  100  &   1.598e-06  &  4.508  &  1.944e-05 & 3.706  &  1.944e-05 & 3.706 \\ 
  200  &   5.794e-08  &  4.786  &  8.602e-07 & 4.498  &  8.603e-07 & 4.498 \\ 
  300  &   7.942e-09  &  4.901  &  1.240e-07 & 4.777  &  1.240e-07 & 4.778 \\ 
  400  &   1.938e-09  &  4.903  &  3.176e-08 & 4.735  &  3.176e-08 & 4.734 \\ 
  500  &   6.658e-10  &  4.788  &  1.105e-08 & 4.732  &  1.105e-08 & 4.732 \\ 
  600  &   2.958e-10  &  4.450  &  4.626e-09 & 4.775  &  4.626e-09 & 4.775 \\ 

		\hline
	\end{tabular}
	\caption{Order of convergence for \eqref{eq:12} with $p=1,\dots,3$ obtained with WENO5 on a cartesian grid with $N\times N$ cells\label{tab:Cos3}}
\end{table}
For the tests run with the $\cos$ vortexes \eqref{eq:12}, we observe that for $p=1$ we reach order 3 for $h$ and around $2.4$ for $u$, a little better than expected, for $p=2$ the order reaches 4 and almost 5 for $h$, and finally for $p=3$ we essentially obtain the expected 5th order of the scheme already at the mesh $300\times 300$.

\begin{table}
	\centering
	p=2\\
	\begin{tabular}{|c||c|c||c|c||c|c|}
		\hline
		  N   & Error h  &  Order h & Error u  &  Order u & Error v  &  Order v \\ 
   25  &   3.367e-04  &  0.000  &  2.965e-03 & 0.000  &  2.966e-03 & 0.000 \\ 
   50  &   1.301e-04  &  1.372  &  1.239e-03 & 1.258  &  1.239e-03 & 1.259 \\ 
  100  &   3.142e-05  &  2.050  &  3.213e-04 & 1.948  &  3.213e-04 & 1.948 \\ 
  200  &   3.005e-06  &  3.387  &  5.399e-05 & 2.573  &  5.399e-05 & 2.573 \\ 
  300  &   5.146e-07  &  4.352  &  1.387e-05 & 3.353  &  1.387e-05 & 3.353 \\ 
  400  &   1.205e-07  &  5.045  &  5.169e-06 & 3.430  &  5.169e-06 & 3.430 \\ 
  500  &   3.715e-08  &  5.275  &  2.542e-06 & 3.181  &  2.542e-06 & 3.181 \\ 
  600  &   1.427e-08  &  5.247  &  1.362e-06 & 3.420  &  1.362e-06 & 3.420 \\ 

		\hline
	\end{tabular}\\
	p=3	\\
	\begin{tabular}{|c||c|c||c|c||c|c|}
		\hline
		  N   & Error h  &  Order h & Error u  &  Order u & Error v  &  Order v \\ 
   25  &   3.186e-04  &  0.000  &  2.459e-03 & 0.000  &  2.459e-03 & 0.000 \\ 
   50  &   1.092e-04  &  1.544  &  1.017e-03 & 1.273  &  1.017e-03 & 1.273 \\ 
  100  &   2.372e-05  &  2.203  &  2.401e-04 & 2.083  &  2.401e-04 & 2.083 \\ 
  200  &   1.815e-06  &  3.708  &  3.212e-05 & 2.902  &  3.212e-05 & 2.902 \\ 
  300  &   2.577e-07  &  4.815  &  7.621e-06 & 3.548  &  7.621e-06 & 3.548 \\ 
  400  &   5.728e-08  &  5.227  &  2.966e-06 & 3.281  &  2.966e-06 & 3.281 \\ 
  500  &   1.810e-08  &  5.162  &  1.339e-06 & 3.565  &  1.339e-06 & 3.565 \\ 
  600  &   7.234e-09  &  5.030  &  6.317e-07 & 4.119  &  6.317e-07 & 4.119 \\ 

		\hline
	\end{tabular}\\
	p=4\\
	\begin{tabular}{|c||c|c||c|c||c|c|}
		\hline
		  N   & Error h  &  Order h & Error u  &  Order u & Error v  &  Order v \\ 
   25  &   3.047e-04  &  0.000  &  2.424e-03 & 0.000  &  2.425e-03 & 0.000 \\ 
   50  &   1.016e-04  &  1.584  &  9.224e-04 & 1.394  &  9.224e-04 & 1.394 \\ 
  100  &   2.232e-05  &  2.187  &  2.205e-04 & 2.065  &  2.205e-04 & 2.065 \\ 
  200  &   1.657e-06  &  3.752  &  2.706e-05 & 3.027  &  2.706e-05 & 3.027 \\ 
  300  &   2.273e-07  &  4.899  &  6.324e-06 & 3.585  &  6.324e-06 & 3.585 \\ 
  400  &   5.102e-08  &  5.194  &  2.485e-06 & 3.247  &  2.485e-06 & 3.247 \\ 
  500  &   1.639e-08  &  5.089  &  1.075e-06 & 3.755  &  1.075e-06 & 3.755 \\ 
  600  &   6.636e-09  &  4.959  &  4.934e-07 & 4.272  &  4.934e-07 & 4.272 \\ 

		\hline
	\end{tabular}\\
	p=5\\
	\begin{tabular}{|c||c|c||c|c||c|c|}
		\hline
		  N   & Error h  &  Order h & Error u  &  Order u & Error v  &  Order v \\ 
   25  &   2.920e-04  &  0.000  &  2.256e-03 & 0.000  &  2.255e-03 & 0.000 \\ 
   50  &   1.002e-04  &  1.543  &  9.002e-04 & 1.325  &  9.002e-04 & 1.325 \\ 
  100  &   2.275e-05  &  2.139  &  2.180e-04 & 2.046  &  2.180e-04 & 2.046 \\ 
  200  &   1.749e-06  &  3.701  &  2.608e-05 & 3.064  &  2.608e-05 & 3.064 \\ 
  300  &   2.413e-07  &  4.886  &  6.022e-06 & 3.615  &  6.021e-06 & 3.615 \\ 
  400  &   5.463e-08  &  5.164  &  2.383e-06 & 3.222  &  2.383e-06 & 3.222 \\ 
  500  &   1.766e-08  &  5.060  &  1.027e-06 & 3.771  &  1.027e-06 & 3.770 \\ 
  600  &   7.171e-09  &  4.944  &  4.691e-07 & 4.300  &  4.691e-07 & 4.300 \\ 

		\hline
	\end{tabular}
	\caption{Order of convergence for \eqref{eq:hvortexCInfp} with $p=2,\dots,5$ obtained with WENO5 on a cartesian grid with $N\times N$ cells\label{tab:Cinf}}
\end{table}

For the $\mathcal C^\infty$ cases we observe that a finer mesh is needed to reach the expected accuracy, in particular for the \eqref{eq:hvortexCInfp} tests, where we barely reach it for $u$ with the mesh with $N=600$. We observe in this case that we obtain better convergence orders for smaller times, e.g. $T=0.1$. 

\begin{table}
	\centering
	p=2\\
	\begin{tabular}{|c||c|c||c|c||c|c|}
		\hline
		  N   & Error h  &  Order h & Error u  &  Order u & Error v  &  Order v \\ 
   25  &   3.315e-04  &  0.000  &  2.935e-03 & 0.000  &  2.935e-03 & 0.000 \\ 
   50  &   1.271e-04  &  1.383  &  1.216e-03 & 1.271  &  1.216e-03 & 1.272 \\ 
  100  &   3.040e-05  &  2.064  &  3.141e-04 & 1.952  &  3.141e-04 & 1.952 \\ 
  200  &   2.895e-06  &  3.392  &  5.286e-05 & 2.571  &  5.286e-05 & 2.571 \\ 
  300  &   4.950e-07  &  4.356  &  1.358e-05 & 3.352  &  1.358e-05 & 3.352 \\ 
  400  &   1.160e-07  &  5.044  &  5.067e-06 & 3.426  &  5.067e-06 & 3.426 \\ 
  500  &   3.579e-08  &  5.269  &  2.492e-06 & 3.180  &  2.492e-06 & 3.180 \\ 
  600  &   1.377e-08  &  5.239  &  1.335e-06 & 3.424  &  1.335e-06 & 3.424 \\ 

		\hline
	\end{tabular}\\
	p=3	\\
	\begin{tabular}{|c||c|c||c|c||c|c|}
		\hline
		  N   & Error h  &  Order h & Error u  &  Order u & Error v  &  Order v \\ 
   25  &   2.897e-04  &  0.000  &  2.259e-03 & 0.000  &  2.259e-03 & 0.000 \\ 
   50  &   8.398e-05  &  1.787  &  8.381e-04 & 1.431  &  8.381e-04 & 1.430 \\ 
  100  &   1.536e-05  &  2.451  &  1.734e-04 & 2.273  &  1.734e-04 & 2.273 \\ 
  200  &   1.040e-06  &  3.884  &  2.194e-05 & 2.982  &  2.194e-05 & 2.982 \\ 
  300  &   1.412e-07  &  4.925  &  5.295e-06 & 3.506  &  5.295e-06 & 3.506 \\ 
  400  &   3.214e-08  &  5.146  &  2.072e-06 & 3.262  &  2.072e-06 & 3.262 \\ 
  500  &   1.036e-08  &  5.074  &  8.860e-07 & 3.806  &  8.860e-07 & 3.806 \\ 
  600  &   4.199e-09  &  4.952  &  4.054e-07 & 4.288  &  4.054e-07 & 4.288 \\ 

		\hline
	\end{tabular}\\
	p=4\\
	\begin{tabular}{|c||c|c||c|c||c|c|}
		\hline
		  N   & Error h  &  Order h & Error u  &  Order u & Error v  &  Order v \\ 
   25  &   2.637e-04  &  0.000  &  2.123e-03 & 0.000  &  2.123e-03 & 0.000 \\ 
   50  &   6.696e-05  &  1.978  &  6.581e-04 & 1.690  &  6.581e-04 & 1.690 \\ 
  100  &   1.052e-05  &  2.670  &  1.210e-04 & 2.444  &  1.210e-04 & 2.444 \\ 
  200  &   5.934e-07  &  4.148  &  1.284e-05 & 3.236  &  1.284e-05 & 3.235 \\ 
  300  &   7.766e-08  &  5.015  &  3.256e-06 & 3.385  &  3.256e-06 & 3.384 \\ 
  400  &   1.812e-08  &  5.059  &  1.156e-06 & 3.601  &  1.156e-06 & 3.601 \\ 
  500  &   6.034e-09  &  4.928  &  4.470e-07 & 4.257  &  4.470e-07 & 4.257 \\ 
  600  &   2.515e-09  &  4.799  &  1.982e-07 & 4.461  &  1.982e-07 & 4.461 \\ 

		\hline
	\end{tabular}\\
	p=5\\
	\begin{tabular}{|c||c|c||c|c||c|c|}
		\hline
		  N   & Error h  &  Order h & Error u  &  Order u & Error v  &  Order v \\ 
   25  &   2.020e-04  &  0.000  &  1.739e-03 & 0.000  &  1.739e-03 & 0.000 \\ 
   50  &   6.518e-05  &  1.632  &  5.841e-04 & 1.574  &  5.841e-04 & 1.574 \\ 
  100  &   8.822e-06  &  2.885  &  9.722e-05 & 2.587  &  9.721e-05 & 2.587 \\ 
  200  &   4.440e-07  &  4.312  &  8.910e-06 & 3.448  &  8.909e-06 & 3.448 \\ 
  300  &   5.833e-08  &  5.006  &  2.335e-06 & 3.302  &  2.335e-06 & 3.302 \\ 
  400  &   1.396e-08  &  4.972  &  7.587e-07 & 3.908  &  7.587e-07 & 3.908 \\ 
  500  &   4.712e-09  &  4.866  &  2.826e-07 & 4.425  &  2.826e-07 & 4.425 \\ 
  600  &   1.972e-09  &  4.776  &  1.250e-07 & 4.476  &  1.250e-07 & 4.476 \\ 

		\hline
	\end{tabular}
	\caption{Order of convergence for \eqref{eq:hvortexCTan} with $p=2,\dots,5$ obtained with WENO5 on a cartesian grid with $N\times N$ cells\label{tab:CTan}}
\end{table}

In the $\arctan$ vortex \eqref{eq:hvortexCTan} we have smaller amplitude of higher derivatives, in particular for large $p$, so with $p=3,\,4$ we see an order of accuracy larger than 4 for $N=600$, which is a bit better than with the only exponential test, still with a quite fine mesh. The same reasoning holds for these tests and with smaller final times we can observe better convergence order (around 4.8 for $p\geq 3$).

\section{Conclusion}
Steady and moving vortex solutions are useful tools to test the order of accuracy of high order methods for hyperbolic balance laws (shallow water and Euler equations). The vortexes known in literature are, anyway, either discontinuous in some derivatives \cite{ricchiuto2009stabilized}, hence not suitable for very high order methods, or non compact supported \cite{shu1998essentially}, hence presenting troubles with boundary conditions if the parameters are not carefully chosen.

In this work we propose a class of $\mathcal{C}^{2p}$ compactly supported vortexes with arbitrary $p$, and some $\mathcal{C}^\infty$ compactly supported vortexes.
The latter are very attractive as they do not need to compute integrals, but only derivatives, and have all the derivatives continuous. On the other side, their derivatives are particularly large and the expected order of accuracy is reached only for very fine meshes. The one that obtain better results in this class is the one based on the $\arctan$ with $p=5$ for final time $T=0.1$. For larger times the order is not so neat.

Among the $\cos$ based vortexes, the one with $p=3$ gives better results (for order 5 schemes) even for coarser meshes $N_x\approx 300$, as its derivative are not too large. On the other side, this vortex can only be used to test methods up to order 6 of accuracy.

Overall, we suggest the use of $\mathcal{C}^\infty$ vortexes for testing arbitrarily high order methods, while, for a fixed method of order $d\leq 2p$ the recipe \eqref{eq:12} builds a vortex where the order of accuracy can be observed in coarser meshes.

\section*{Acknowledgments}
D.T. acknowledges Wasilij Barsukow for a fruitful discussion on vortexes and their origin.

\appendix
\section{Values of $H_p$}\label{app:exactCos}
Here, we provide the explicit form of $H_p$ for $p=1,\dots, 3$, the corresponding $h$ can be computed from \eqref{eq:12}. The value of $\Gamma$ can be set given $h_{min}$, using 
\begin{equation}
	\Gamma = \frac{\pi}{2^p r_0}\sqrt{\frac{g(h_0-h_{min})}{H_p(\pi/2)-H_p(0)}}.
\end{equation}
We provide in the following also the value of $H_p(\pi/2)-H_p(0)$.

\begin{align*}
	H_1(x)=&\frac{\cos(2x)}{8} + \frac{x\sin(2x)}{4} + \frac{\cos(2x)^2}{64} + \frac{3x^2}{16} + \frac{x\cos(2x)\sin(2x)}{16},\\
	H_1(\pi/2)-&H_1(0) =  \frac{3\pi^2}{64}- \frac{1}{4}  ,\end{align*}
\begin{align*}
	H_2(x) =& \frac{35\cos(2x)}{384} + \frac{35x\sin(2x)}{192 }+ \cos(x)^6\left( \frac{\cos(x)^2}{64 }+ \frac{7}{288}\right) + \frac{35\cos(2x)^2}{3072 }+ \\ &\frac{35x^2}{256 }+ \frac{35x\cos(2x)\sin(2x)}{768 } + \frac{x\cos(x)^5\sin(x)(\cos(x)^2 + 7/6)}{8},\\
	H_2(\pi/2)-&H_2(0) = \frac{35\pi^2}{1024 }- \frac{2}{9}  ,\end{align*}
\begin{align*}
	H_3(x)=&	\frac{77\cos(2x)}{1024 }+ \frac{77x\sin(2x)}{512 }+ \frac{33\cos(x)^6(\cos(x)^2/64 + 7/288)}{40} + \cos(x)^{10}\left(\frac{\cos(x)^2}{144} + \frac{11}{1200}\right) +\\
	 &\frac{77\cos(2x)^2}{8192} + \frac{231x^2}{2048 } + \frac{77x\cos(2x)\sin(2x)}{2048 }+ \frac{33x\cos(x)^5\sin(x)(\cos(x)^2 + 7/6)}{320 }+\\ & \frac{x\cos(x)^9\sin(x)(\cos(x)^2 + 11/10)}{12},\\
	 H_3(\pi/2)-&H_3(0) =   \frac{231\pi^2}{8192 }- \frac{359}{1800  }.
\end{align*}

\section{Runge Kutta (6,5)}\label{sec:RK65}
The Butcher Tableau of Butcher's Runge Kutta (6,5) method  used is the following
\begin{align}
	\renewcommand\arraystretch{1.2}
	\begin{array}{c|cccccc}
		0&  &  &  &  &  & \\
		\frac{1}{4} & \frac{1}{4} &  &  &  &  & \\
		\frac{1}{4} & \frac{1}{8} & \frac{1}{8} &  &  &  & \\
		\frac{1}{2} & 0 & 0 & \frac{1}{2} &  &  & \\
		\frac{3}{4} & \frac{3}{16} & - \frac{3}{8} & \frac{3}{8} & \frac{9}{16} &  & \\
		1 & - \frac{3}{7} & \frac{8}{7} & \frac{6}{7} & - \frac{12}{7} & \frac{8}{7} & \\
		\hline
		& \frac{7}{90} & 0 & \frac{16}{45} & \frac{2}{15} & \frac{16}{45} & \frac{7}{90}\\
	\end{array}.
\end{align}
\bibliographystyle{siam}
\bibliography{biblio}

 \end{document}